\magnification=1200
\input amssym.def
\input amssym.tex
 \font\newrm =cmr10 at 24pt
\def\bul{\raise .9pt\hbox{\newrm .\kern-.105em } }

\overfullrule=0pt
\font\ninerm=cmr10 at 9pt
\font\authorfont=cmr10

 \def\fr{\frak}

 \baselineskip 16pt
 
 \def\h{\hbox{ }}
 
 \def\p{{\fr p}}
 \def\u{{\fr u}}
 \def\r{{\fr r}}
 
 \def\m{{\fr m}}
 \def\n{{\fr n}}
 \def\a{{\fr a}}

 \def\ss{{\fr s}}
 \def\k{{\fr k}}
 \def\b{{\fr b}}
 
 \def\hh{{\fr h}}

 \def\g{{\fr g}}

 \def\<{\le}
 \def\>{\ge}

 \def\s{{\h\subset\h}}
 
 \def\vs{\vskip }

 \def\mapright#1
  {\smash{\mathop
  {\longrightarrow}
  \limits^{#1}}}

 \def\kk#1{{\kern .4 em} #1}   
 \def\vs{\vskip 1pc}

\rm 
	\centerline{\bf  On the Centralizer of  $K$ in $U(\g)$}\vskip 1pc

\centerline{\authorfont BERTRAM KOSTANT\footnote*{\ninerm \baselineskip=11pt
Research supported in part by the KG\&G Foundation.}}
 \vskip 1.5pc

\hfill{\it Dedicated  with respect to  Ernest Vinberg}

\hfill{\it  on the occasion of his seventieth birthday}

\vskip 1pc
\noindent {\bf Abstract.} Let $\g = \k +\p$ be a complexified Cartan decomposition
of a complex semisimple Lie algebra $\g$ and let $K$ be the
subgroup of the adjoint group of $\g$ corresponding to $\k$. If 
$H$ is an irreducible Harish-Chandra module of $U(\g)$, then $H$ is
completely determined by the finite-dimensional action of the
centralizer $U(\g)^K$ on any one fixed primary $\k$ component in
$H$. This original approach of Harish-Chandra to a determination of
all $H$ has largely been abandoned because one knows very little
about generators of $U(\g)^K$. Generators of $U(\g)^K$ may be given by
generators of the symmetric algebra analogue $S(\g)^K$.  Let
$S_m(\g)^K,\,\,m\in \Bbb Z_+$, be the subalgebra of $S(\g)^K$
defined by $K$-invariant polynomials of degree at most $m$. For convenience write $A = S(\g)^K$ and
$A_m$ for the subalgebra of $A$ generated by $S_m(\g)^K$. Let
$Q$ and $Q_m$ be the respective quotient fields of $A$ and
$A_m$. We prove that if $n= dim\,\g$ one has $Q= Q_{2n}$. 

We also determine the variety, $Nil_K$, of
unstable points with respect to the action $K$ on $\g$ and show that $Nil_K$ is already
defined by $A_{2n}$. As pointed out to us by Hanspeter Kraft this fact together with a result 
of Harm Derksen (See [D]) implies, indeed, that $A= A_r$ where $r = {2n\choose 2}\,dim\,\p$.

\vskip 1pc
\centerline{\bf 1. Introduction}\vskip 1pc {\bf 1.1.} Let
$\g$ be a complex semisimple Lie algebra. The value taken by the
Killing form,
$B$, on
$w,z\in \g$ will be denoted by 
$(w,z)$. Let
$$\g=
\k +
\p
\eqno (1.1)$$ be a complexified Cartan decomposition and let
$\theta$ be the corresponding complexified Cartan involution. One
has that
$[\p,\p]$ is an ideal of $\k$ (and $\p + [\p,\p]$ is an ideal of
$\g$).
We will assume that (1.1) is proper in the sense that $$\k =
[\p,\p]\eqno (1.2) $$ (i.e., (1.1) arises from the Cartan decomposition of a real
form of $\g$ without ``compact components"). Let
$G$ be the adjoint group of $\g= Lie\,\g$ and let $K\s G$ be the subgroup
corresponding to $\k$. Of course $G$ has trivial center.  

 We recall that the centralizer
$U(\g)^K$ of $K$ in the universal enveloping algebra $U(\g)$ of $\g$ played a key
role in Harish-Chandra's original approach to the study of certain infinite
dimensional representations of 
$\g$. A critical end product of the theory is the existence of irreducible 
Harish-Chandra modules. Such a module $M$ is an irreducible $U(\g)$-module which
not only is completely reducible as a $\k$-module but also the primary
components are finite dimensional. Any such primary component then defines a
finite-dimensional $U(\g)^K$-module and, remarkably, the entire
$U(\g)$-module $M$ is completely determined by the action of $U(\g)^K$ on any one
fixed primary component. An early consequence of all of this is Harish-Chandra's
subquotient theorem. (For a considerable simplification and clarification of
Harish-Chandra's proof see Lepowsky [L] and Lepowsky-McCollum [L-M]. See also
Wallach [W] and Vogan [V-1]). With the the determination of Harish-Chandra
modules reduced to a determination of the finite-dimensional representation
theory of
$U(\g)^K$ one might have expected a subsequent development of representation
theory along these lines. However this has not been the case although a
considerable effort in this direction is seen in [V-1]. The main result of [V-1] is
a classification theorem.  One major obstacle to making progress with this
approach is that the algebra $U(\g)^K$ is poorly understood. This is more or less
attested to by Vogan in [V-2] where he remarks that $U(\g)^K$ is ``hideously
complicated". See p.~17 in [V-2]. Also see [K-T] for a glimpse into this complication.

It is not difficult to construct a linear
basis of $U(\g)^K$. The difficulty lies with its ring structure. Progress would
be made if we could pin down a set of (algebra) generators of $U(\g)^K$. Indeed
focusing on the primary component, given by Vogan's minimal $\k$-type, 
the corresponding representation of $U(\g)^K$ is given by a one-dimensional
character. Consequently the whole $U(\g)$-module $M$ is known as soon as one knows
the scalar values assigned to these generators by the character. 

The
algebra  $U(\g)^K$ has a natural filtration and
$PBW$ implies an algebra isomorphism $$Gr\,U(\g)^K \cong S(\g)^K\eqno (1.3)$$
where $S(\g)^K$ is the finitely generated integral domain of
$Ad\,K$ invariants in the symmetric algebra $S(\g)$. A set of
homogeneous generators of $S(\g)^K$ then yields a set of generators of $U(\g)^K$.
The main results of this paper together with a result of Derksen in [D] yields generators of
$S(\g)^K$.\vs  {\bf 1.2.} The adjoint action of $k\in K$ on $z\in \g$ will be denoted by
$k\cdot z$. If $z\in \g$, then $K\cdot z $ is Zariski closed if
and only if it is closed in the usual Hausdorff topology. Let
$$Cl\,\g =
\{z\in
\g\mid K\cdot z\,\,\hbox{ is closed}\}$$ For notational
convenience we will put
$A = S(\g)^K$. Also for notational convenience we identify $S(\g)$ with the
algebra of polynomial functions on $\g$ where for any $x,y\in \g$ and $m\in \Bbb
Z_+$, one has $x^m(y) = (x,y)^m$. Then
$A$ is the affine algebra of the affine variety $V$ of
all homomorphisms $A\to \Bbb C$, i.e., all closed points in
$Spec\,A$. Then, from invariant theory, one knows that $$V\cong
Cl\,\g/K$$ i.e., $$V\,\,\hbox{identifies with the set of all
closed $K$-orbits in $\g$.}\eqno (1.4)$$\vskip .5pc For any $z\in
Cl\,\g$ we will let $$v_z\in V\,\,\hbox{be the point
corresponding to $K\cdot z.$}\eqno (1.5)$$ \vskip .5pc The
symmetric algebra $S(\g)$ is filtered by the subspaces
$S_m(\g),\,\,m\in
\Bbb Z_+$, where $S_m(\g) = \sum_{j=0}^{m}\,\,S^j(\g)$. Obviously
$$S_m(\g)^K = \sum_{j=0}^{m}\,\,S^j(\g)^K\eqno (1.6)$$ But then 
$A$ is filtered by the subalgebras $A_m,\,m\in \Bbb Z_+$, where we
let
$$A_m\,\,\,\hbox{be the subalgebra of $A$ generated by
$S_m(\g)^K$}\eqno (1.7)$$ Let $V_m$ be the affine variety
corresponding to $A_m$. The injection
$$0\longrightarrow A_m\longrightarrow A\eqno (1.8)$$ defines a
dominant morphism $$\gamma_{m}:V\to V_m\eqno (1.9)$$\vskip .5pc Let
$Q\,\,(\hbox{resp}. Q_m)$ be the quotient field of
$A\,\,(\hbox{resp}. A_m)$ and let $$n= dim\,\g \eqno (1.10)$$ 

The first main result is \vs {\bf Theorem 1.1}. {\it The dominant morphism
$\gamma_{2n}$ is birational so that $$Q = Q_{2n}\eqno (1.11)$$ In
particular any $h\in A$ is of the form $$h = f/g\eqno (1.12)$$
where $f,g\in A_{2n}$ and of course $g\neq 0$.}\vs {\bf 1.3.}  Let $z\in \g$ be
arbitrary. Then $z$ can be uniquely written $$z = x +
y,\,\,\hbox{where}\,\,x\in \k\,\,\hbox{and}\,\,y\in \p\eqno (1.13)$$
 Let $\g(z)$ be the Lie
subalgebra of $\g$ generated by $x$ and $y$. We will use this notation throughout the
paper.

In constrast to the closed $K$-orbits in $\g$, consider the cone of
$K$-unstable points in $\g$. Let
$$Nil_K = \{z\in \g \mid f(z)=0,\,\,\forall\,\, \hbox{homogeneous $f\in S(\g)^K$
of positive degree}\}$$ Since $S(\g)^G\s S(\g)^K$ obviously $Nil_K$ is a
subvariety of the nilcone of $\g$. \vs {\bf Theorem 1.2.}
{\it Let
$z\in
\g$. Then $z\in Nil_K$ if and only if $\g(z)$ is a (nilpotent) Lie algebra of
nilpotent elements.}\vs For a number of results about the 
nilcones of the actions of $K$, or rather $K_{\theta}$, (defined in (2.32) below)
on multiple copies of $\p$ see [K-W]. Also see [P-3]. For the
case we are considering here, Wallach raised the
question for a determination of some value of
$m\in
\Bbb Z_+$ with the property that
$Nil_K$ is given already by the homogeneous elements in $A_m$ of positive degree. The
following result answers this question with the same value of $m$ appearing in Theorem
1.1, namely
$m=2n$. \vs {\bf Theorem 1.3.} {\it Let
$z\in
\g$. Then
$z\in Nil_K$ if and only if
$$f(z) = 0,\,\,\,\forall f\in A_{2n}\,\,\hbox{of positive degree.}$$} 
\vskip .5pc The idea of using a degree which defines $Nil_K$ (in this case $2n$) to determine $r$
such that $A= A_r$, goes back to Popov. See [P-1] and [P-2]. Harm Derksen in [D] has sharply reduced
Popov's estimate of $r$ . Thus combining Theorem 1.3 with the result in [D] one has\vs {\bf Theorem
1.4} {\it One has $$A= A_r\eqno (1.14)$$ where $$r = {2n\choose 2}\,dim\,\p\eqno (1.15)$$ where, we
recall $n=dim\,\g$. }\vs I thank Hanspeter Kraft for informing me about Derksen's result. Kraft
formulated Theorem 1.4, seeing it as an immediate consequence of my Theorem 1.3 and
Derksen's result. I also thank Nolan Wallach for motivating me to think
about finding an integer $m$ such that $A_m$ defines $Nil_K$ (see Theorem 1.3). I also thank him for many
conversations about invariant theory.
\vskip 1.5pc
\centerline{\bf 2. The proof of Theorems 1.1, 1.2, 1.3 and 1.4}\vskip 1pc
{\bf 2.1.} Let $\Phi = \Phi(X,Y)$ be the free Lie algebra, over
$\Bbb C$ on two generators  $X,Y$. The Lie algebra $\Phi$ is
naturally graded over $\Bbb Z_+$ with homogeneous spaces $\Phi^j$.
It is then clearly filtered by the subspaces $\Phi_m,\,\,m\in \Bbb
Z_+$, where
$$\Phi_{m} = \sum_{j=0}^m\,\Phi^j\eqno (2.1)$$ Clearly $$\Phi_{m+1} = \Phi_m +
[X,\Phi_m] + [Y,\Phi_m]\eqno (2.2)$$\vskip .5pc Using notation introduced in \S 1.3
one then has a Lie algebra epimorphism, $$\xi_z:\Phi\to
\g(z),\,\,\,\hbox{where}\,\,\xi_z(X) =
x\,\,\,\hbox{and}\,\,\xi_z(Y)=y\eqno (2.3)$$\vskip .5pc The Lie
subalgebra $\g(z)$ of $\g$ is filtered by the subspaces $\g_m(z)$
where we put $\g_m(z) = \xi_z(\Phi_m)$. By (2.2) one has
$$\g_{m+1}(z) = \g_m(z) + [x,\g_m(z)] + [y,\g_m(z)]\eqno
(2.4)$$\vskip .5pc {\bf Proposition 2.1.} {\it For any $z\in \g$ one
has
$$\g_{n-1}(z) = \g(z)\eqno (2.5)$$}\vs {\bf Proof.} It follows
immediately from (2.4) that $\g(z) = \g_{m}(z)$ in case $$\g_{m}(z)
= \g_{m+1}(z)\eqno (2.6)$$ Indeed (2.6) implies that $\g_k(z) =
\g_m(z)$ for all $k\in \Bbb Z_+$ where $k\geq m$. 

The statement of the proposition is obviously true if
$dim\,\g_1(z)\leq 1$. We can therefore assume $dim\,\g_1(z) = 2$.
We refer to the equality (2.6) as ``stability at $m$". If one does
not have stability at $m$ then clearly $$dim\,\g_{m+1}(z) > m+1\eqno
(2.7)$$ But then nonstability at $n-1$ yields the contradictory
statement $dim\,\g_{n}(z)> n = dim\,\g$. Hence one necessarily
has stability at $n-1$. QED\vs {\bf 2.2.} 
If $z= x+y$ is the decomposition (1.13) for $z\in\g$, then
obviously $k\cdot z  = k\cdot x + k\cdot y$ is the decomposition
(1.13) for $k\cdot z$ for any $k\in K$. The following simple
statement is important for us.
\vs {\bf Proposition 2.2.} {\it Let $T,T'\in \Phi_n$. Then
$f_{T,T'}\in S_{2n}(\g)^K$ where, for $z\in \g$, $$f_{T,T'}(z) =
(\xi_z(T),\xi_z(T'))\eqno (2.8)$$}\vs {\bf Proof.} We only
have to observe that $f_{T,T'}\in S_{2n}(\g)$. The remainder
follows from invariance of the Killing form and the fact that for
$W\in \Phi,\,z\in \g$ and $k\in K$,$$k\cdot\xi_z(W) =
\xi_{k\cdot z}(W)\eqno (2.9)$$ QED\vs Let $$\g^{K\,reg}=\{z\in
\g\mid \g(z) = \g\}$$ Thus, by Proposition 2.1, $z\in \g^{K\,reg}$
if and only if
$$\eqalign{\g_{n-1}(z) &= \g_{n}(z)\cr &=\g\cr}\eqno (2.10)$$ One
readily constructs some $z\in \g$ to show that $\g^{K\,reg}$ is not
empty. See Appendix for a proof that  $\g^{K\,reg}$ is not empty.

 Let $d(n) =
dim\,\Phi_n$. Let $T_j,\,j=1,\ldots,d(n)$, be a basis of $\Phi_n$.
The following is a restatement of Proposition 2.1 and (2.10). \vs
{\bf Proposition 2.3.} {\it Let $z\in \g$. Then
$\xi_z(T_j),\,j=1,\ldots,d(n)$, spans $\g(z)$. In particular
$z\in
\g^{K\,reg}$ if and only if
$\xi_z(T_j),\,j=1,\ldots, d(n)$, spans $\g$.} \vs As
functions on $\g$ the entries of
the 
$d(n)\times d(n)$ matrix $M(z)$ given by $$M_{i\,j}(z) =
(\xi_z(T_i),\xi_z(T_j))$$ are in
$S_{2n}(\g)^K$. 

For any $z\in \g$ let $K_z$ be the stabilizer of $z$ with respect
to the adjoint action $K$ on $\g$ . Let $\k_z = Lie\, K_z$.
Clearly $$\k_z\,\,\hbox{is the centralizer of $\g(z)$ in $\k$}\eqno
(2.11)$$ From the semisimplicity of $\g$ one then has $$\k_z =
0\,\,\hbox{for any $z\in \g^{K\,reg}$}\eqno (2.12)$$\vs {\bf Theorem
2.4.} {\it $\g^{K\,reg}$ is a nonempty Zariski open subset of $\g$. Furthermore if
$z\in
\g^{K\,reg}$ then the $K$-orbit $K\cdot z$ is closed. That is,
$$\g^{K\,reg}\s Cl(\g)\eqno (2.13)$$ Put $$V^{K\,reg} =\{v\in V\mid
v=v_z\,\,\hbox{for some $z\in \g^{K\,reg}$}\}\eqno (2.14)$$ Then
$V^{K\,reg}$  is a nonempty Zariski open (and hence dense) subset
of $V$.} \vs {\bf Proof.} Let $z\in
\g$. Then clearly
$$rank\,M(z)\leq\dim\,g(z)\eqno (2.15)$$ But since the Killing form
is nonsingular on $\g$ it follows that $$rank\,M(z) =
dim\,\g\,\,\iff z\in \g^{K\,reg}$$ 

Let $z\in \g^{K\,reg}$ and let $z'\in \overline{K\cdot z}$. But
then clearly $M(z) = M(z')$ so that $z'\in \g^{K\,reg}$. But then
$k_{z'}= 0$ by (2.12). Thus $dim\,K\cdot z = dim\,K\cdot z'$. This
implies that $K\cdot z$ is closed since the $K$-orbits on the
boundary of $K\cdot z$ must have dimension smaller than
$dim\,K\cdot z$. But now the determinants of all the
$dim\,\g\times dim\,\g$ minors of $M(z)$ are in $A$. It is an easy exercise 
to show that $\g^{K\,reg}$ is not
empty. (As mentioned above a proof that $\g^{K\,reg}$ is not
empty is given in the Appendix.) This proves that $\g^{K\,reg}$
is a nonempty Zariski open subset of $\g$ and $V^{K\,reg}$ is a nonempty Zariski
open subset of $V$. QED\vs {\bf Remark 2.5.} Note that since the entries of
$M(z)$ are in
$S_{2n}(\g)^K$ the determinants of all the $dim\,\g\times dim\,\g$ minors of
$M(z)$ are, in fact, in
$A_{2n}$.\vs {\bf 2.3.} {\bf Proof of Theorem
1.1.}  To show that $\gamma_{2n}$ is birational it suffices, 
by Theorem 2.4, to prove that there exists a nonempty Zariski open subset $V_*\s V$
such that the restriction
$$\gamma_{2n}:V_*
\to V_{2n}\eqno (2.16)$$ is injective. Theorem 2.4 asserts that $V^{K\,reg}$ is
a nonempty open subvariety of $V$. The variety $V_*$, to be constructed, will in
fact be a nonempty open subvariety of $V^{K\,reg}$. Before constructing $V_*$
we will first establish certain properties of the restriction
$$\gamma_{2n}:V^{K\,reg}
\to V_{2n}\eqno (2.17)$$ 
\vskip .5pc Let
$z,z'\in
\g^{K\,reg}$ be such that
$$f(z) = f(z'),\,\,\,\forall\,f \in A_{2n}\eqno (2.18)$$ We will prove that
there exists an automorphism $\pi$ of $\g$, which commutes with $\theta$
such that $z'= \pi(z)$. 

Assume (2.18) is satisfied.
For
$T\in
\Phi_m$ and $j=1,\ldots,d(n)$, let $f_{T,j}\in A_{2n}$ be defined by
putting, for any $w\in \g$, $$f_{T,j}(w) =
(\xi_w(T),\xi_w(T_j))\eqno (2.19)$$ But since $f_{T,j}\in
A_{2n}$, one has $$f_{T,j}(z) = f_{T,j}(z')\eqno (2.20)$$ We
construct a linear isomorphism $$\pi:\g\to \g\eqno (2.21)$$ as
follows: Let $w\in \g$. Then, by (2.10), there exists $T\in \Phi_m$
(obviously not necessarily unique) such that  $\xi_z(T) = w$.
Define (to be shown to be well-defined) $$\pi(w) =
w',\,\,\,\,\hbox{where}\,\,\,w' = \xi_{z'}(T)\eqno (2.22)$$ To
see that $\pi$ is well-defined we have only to establish that if
$T\in \Phi_m$, then $$\xi_{z}(T)=
0,\,\,\,\iff\,\,\xi_{z'}(T)=0\eqno (2.23)$$ But one has
$$\xi_{z}(T)= 0,\,\,\,\iff\,\,\,f_{T,j}(z)
=0,\,\,\,\forall\,\,j=1,\ldots,d(n)\eqno (2.24)$$ The same statement holds when $z'$ replaces $z$. But
then one has (2.23) so that the linear isomorphism $\pi$ is well-defined, noting also that
$$\pi(z) = z'\eqno (2.25)$$\vskip .5pc {\bf Lemma 2.6.} {\it $\pi$
is a Lie algebra automorphism which also commutes with $\theta$. That is,
$\pi$ stabilizes both $\k$ and $\p$.} \vs {\bf Proof.} Let $$\u =
\{t\in \g\mid \pi([t,w]) = [\pi(t),\pi(w)],\,\,\,\forall w\in \g\}$$
Then the Jacobi identity immediately implies that $\u$ is a Lie
subalgebra of $\g$.  Let $w\in \g$ be arbitrary. By (2.10) there
exists $T\in \Phi_{n-1}$ such that $\xi_z(T) = w$. Let $T_X=
[X,T]$ so that $T_X\in \Phi_n$. Define $T_Y\in \Phi_n$ similarly
where $Y$ replaces $X$. Then
$$\eqalign{\xi_z(T_X)&=[x,w]\cr\xi_z(T_Y)&= [y,w]\cr}$$ Let
$\xi_{z'}(T) = w'$ so that $\pi(w) = w'$. Also let $z'=x'+y'$ be
the decomposition (1.13) when $z'$ replaces $z$. Then
$$\eqalign{\xi_{z'}(T_X)&=[x',w']\cr\xi_{z'}(T_Y)&=
[y',w']\cr}$$ Thus the Lie subalgebra $\u$ of $\g$ contains
$x$ and $y$. But then $\u = \g$ since $x$ and $y$ generate $\g$.
Hence $\pi$ is an automorphism. Now let $m\leq n$ where $m\in \Bbb
Z_+$. Let $t_i\in\g,\,i=1,\ldots,m$, where $t_i\in \{x,y\}$. Let
$$w= [t_1,[t_2,[\cdots [t_{m-1},t_m]\cdots ]$$ Then note that
$w\in\k$ or $\p$ according as the number indices
$j$ such that $t_j=y$ is even or odd. It follows immediately that
$\pi$ stabilizes both $\k$ and $\p$. QED\vs We will next restrict
$\gamma_{2n}$ to a nonempty Zariski open subset $V_1$ of $V^{K\,reg}$ to guarantee
that
$\pi$ is an inner automorphism. 

One knows the degrees of the generators of $S(\g)^G$. The maximum
degree is the Coxeter number of some simple component of $\g$.
This number is certainly less than $n$ and hence $$S(\g)^G\s
A_{2n}\eqno (2.26)$$ 

Let $\Gamma$ be the quotient of the group $Out\,\g$ of outer
automorphisms of $\g$ by the normal subgroup $Inn\,g = G$ of inner
automorphisms. The group
$\Gamma$ is finite. The image, in $\Gamma$, of any $\alpha\in Out\,g$ will be
denoted by $\sigma_{\alpha}$. Clearly $S(\g)^G$ is stable under the action of
$Out\,G$ on $S(\g)$. But this clearly defines a representation of
$$\Gamma \to Aut\, S(\g)^G\eqno (2.27)$$ The following is well known but we
will give a proof for completeness.\vs {\bf Lemma 2.7.} {\it The representation
(2.27) is faithful.} \vs {\bf Proof.} Let $\alpha\in Out\,\g$ and assume that
$\alpha\,\notin Inn\,\g$. Let $g\in G$ and put $\alpha' = Ad\,g\,\,\circ\,\,\alpha$.
Then $\sigma_{\alpha}= \sigma_{\alpha'}\neq 1$. However $g$ can be chosen
so that $\alpha'$ stabilizes the Weyl chamber $C$ of a split Cartan subalgebra of
a split real form of $\g$ and $\alpha'|C$ does not reduce to the identity.
However from Weyl group theory one knows that $S(\g)^G$ separates the points
of $C$. This proves that the image of $\sigma_{\alpha}$ in (2.27) is not the
identity. QED \vs For any $1\neq \sigma\in \Gamma$ choose $f_{\sigma}\in S(\g)^G$
such that $f\neq f_{\sigma}$ and let $$F = \prod_{\sigma\in
\Gamma/\{1\}}\,\,(f_{\sigma}-\sigma(f_{\sigma}))\eqno (2.28)$$ putting  $F =1$ if
$\Gamma$ reduces to the identity. Obviously $F\in S(\g)^G\s A_{2n}$. Let
$$\g^{K\,reg}_1 = \{z\in \g^{K\,reg}\mid F(z)\neq 0\}\eqno (2.29)$$ so that
$\g^{K\,reg}_1$, by Theorem 2.4, is a nonempty Zariski open
subset of $\g^{K\,reg}$ and
$$V^{K\,reg}_1 =\{v\mid v= v_z\,\,\hbox{for some}\,\,\,z\in \g^{K\,reg}_1\}\eqno
(2.30)$$ is a nonempty Zariski open subset of
$V^{K\,reg}$. Here we are implicitly using the fact that the intersection of
two nonempty Zariski open subsets of an irreducible variety is again a nonempty
Zariski open set. \vs {\bf Lemma 2.8.} {\it Let
$z,z'\in
\g^{K\,reg}_1$ and assume that (2.18) is satisfied. Let $\pi$ be the
$\g$-automorphism of Lemma 2.6. Then
$\pi$ is inner. That is, $\pi = Ad\,g$ for some $g\in G$ such that $Ad\,g$
stabilizes both $\k$ and $\g$.} \vs {\bf Proof.} If $\pi$ is inner there is
nothing to prove. Assume $\pi$ is not inner and let $1\neq \sigma\in\Gamma$ be
defined by putting $\sigma = \sigma_{\pi^{-1}}$. But by (2.25) one has
$$f_{\sigma}(\pi(z)) = f_{\sigma}(z)\eqno (2.31)$$ But
$$\eqalign{f_{\sigma}(\pi(z)) &= (\pi^{-1} f_{\sigma})(z)\cr &=
(\sigma\,f_{\sigma})(z)\cr}$$ But $(\sigma\,f_{\sigma})(z)\neq f_{\sigma}(z)$
since $F(z)\neq 0$. This contradicts (2.31). Thus $\pi$ is inner. QED\vs Let the
notation be as in Lemma 2.8. We will now restrict $\gamma_{2n}$ even further to
finally guarantee that $g\in K$. 

Taking notation from [K-R] let $$K_{\theta} = \{g\in G\mid Ad\,g\,\,\hbox{stabilizes
both
$\k$ and
$\p$}\}\eqno (2.32) $$ so that, in the notation of Lemma 2.8, $g\in K_{\theta}$. 
 Obviously $K\s K_{\theta}$. Let $Out_{G}\k$ be the group of all
automorphisms of $\k$ of the form $Ad\,g|\k$ for $g\in K_{\theta}$ and let $Inn\,\k$
be the group of all inner automorphisms of $\k$. Obviously $Inn\,\k$ is a
normal subgroup of $Out_{G}\k$. One knows that the quotient group $\Gamma_K =
Out_{G}\k/Inn\,\k$ is finite. See Proposition 1, p. 761 in [K-R]. The
argument yielding (2.26) readily also implies $$S(\k)^K\s A_{2n}\eqno (2.33)$$ Also
the natural action of $Out_{G}\k$ on $S(\k)^K$ descends to a representation
$$\Gamma_K\to Aut\,S(\k)^K\eqno (2.34)$$ The argument establishing Lemma 2.7 is
readily modified (to deal with the case where $\k$ is only reductive but not
semisimple) so that one has \vs {\bf Lemma 2.9.} {\it The representation
(2.34) is faithful.} \vs For each $1\neq \tau\in \Gamma_K$ let $f_{\tau}\in
S(\k)^K$ be such that $f_{\tau}\neq \tau\,f_{\tau}$. If $\Gamma_K$ reduces to the
identity put $F_K = 1$, otherwise let $$F_K =\prod_{\tau\in
\Gamma_K/\{1\}}\,(f_{\tau}-\tau\,f_{\tau})\eqno (2.35)$$\vskip .5pc Let
$$\g_*=
\{z\in \g^{K\,reg}_1\mid F_K(z)\neq 0\}\eqno (2.36)$$ and let $$V_*=\{v\in V\mid v
= v_z\,\,\hbox{for some}\,\,z\in \g_*\}\eqno (2.37)$$ Again, since  
 the intersection of two nonempty Zariski open
subsets of an irreducible variety is again a nonempty Zariski open set, it follows
that $\g_* $ is a nonempty Zariski open subset of
$\g$ and $V_*$ is a nonempty Zariski open subset of
$V$. The following lemma establishes Theorem 1.1. \vs {\bf
Lemma 2.10.} {\it Let
$z,z'\in
\g^*$ be such that $$f(z) = f(z')\eqno (2.38)$$ for all $f\in A_{2n}$. Let $g\in
G$ be given by Lemma 2.8 so that $$Ad\,\, g (z) = z'\eqno (2.39)$$ and $g\in
K_{\theta}$ using the notation of (2.32). Then $g\in K$ so that $$z'\in K\cdot z
\eqno (2.40)$$ proving the injectivity of (2.16) and as, noted in the
beginning of \S 2.3, proving Theorem 1.1.}\vs {\bf Proof.} We first prove that
$Ad\,g|\k\in Inn\,\k$. Assume this is not the case and let
$1\neq \tau$ be the image of $Ad\,g^{-1}|\k$ in $\Gamma_K$. Then, by (2.38),
$$f_{\tau}(Ad\,g\,(z)) = f_{\tau}(z)\eqno (2.41)$$ But, recalling (2.2),
$$\eqalign{ f_{\tau}(Ad\,g\,(z))&= f_{\tau}(Ad\,g\,(x))\cr
&= (Ad\,g^{-1}\,f_{\tau})(x)\cr &= (\tau\,f_{\tau})(x)\cr &= 
(\tau\,f_{\tau})(z)\cr}$$ But this contradicts (2.41) since $F_K(z)\neq 0$.
Hence there exists $k\in K$ such that if $b= k^{-1}\,g$, then $b$ centralizes $\k$.
 But then both the semisimple element
$b_s$ and the unipotent element $b_u$ centralize $\k$ where
$b=b_s\,b_u$ is the Jordan decomposition of $b$. But, as one knows,
the centralizer of $\k$ in $\g$ is commutative, reductive and
contained in $\k$. This readily implies that $b_u =1$ since the nilpotent element
$log\,\,b_u$ must commute with $\k$. Thus
$b$ is semisimple. Hence $b$ centralizes a Cartan subalgebra $\hh$ of
$\g$. Let $\g^b$ be the centralizer of $b$ in $\g$ so that $\hh +\k\s
\g^b$. For any simple component $\g_i$ of $\g$ let $\k_i = \g_i\cap
\g^b$ and let $\p_i$ be the Killing form orthocomplement of $\k_i$
in $\g_i$. Since $\g^b$ contains $\hh$ it is clear that $\g^b$ is
the sum of its intersections with all the simple components of
$\g$. It follows then that $\p_i$ is Killing form orthogonal to
$\k$ so that $\p_i\s \p$. Hence $\p_i + [\p_i,\p_i]$ is an ideal in
$\g_i$. By simplicity either $\p_i = 0$ in which case $\g_i =
\k_i$ so that $\g_i$ makes no nontrivial contribution to $b$ or
$[\p_i,\p_i] = \k_i\s \k$. Since $b$ is in the subgroup of $G$
corresponding to $\hh$ it is then clear that $b\in K$ and hence
$g\in K$. QED 
\vs {\bf 2.4.} {\bf Proof of Theorem 1.2.}
Let $z\in \g$. Then one knows from invariant theory that $K\cdot z$ has a unique
closed $K$-orbit in its  closure (this is immediate from (1.4)). Consequently
$z\in Nil_K$ if and only if $$0\in \overline {K\cdot z}\eqno (2.42)$$ Assume that
$z\in Nil_K$ and let $k_m\in K,\,m\in \Bbb Z_+$, be a sequence such that
$k_m\cdot z$ converges to $0$. But then recalling the decomposition (1.13) one
must have that both $k_m\cdot x$ and $k_m\cdot y$ also converge to $0$. But then
obviously $k_m\cdot w$ converges to $0$ for any $w\in \g(z)$. But then
(recalling that $S(\g)^G\s S(\g)^K$) $w$ is nilpotent for any $w \in \g(z)$. 

Conversely, assume that every element in $\g(z)$ is nilpotent. Then there exists a
Borel subalgebra $\b$ of $\g$ such that $\g(z)\s \n$ where $\n$ is the nilradical
of $\b$. Put $\b' = \theta (\b)$ so that $\theta(\n) = \n'$ where $\n'$ is the
nilradical of $\b'$. Let $\ss = \b\cap \b'$ so that $\ss$ is a
solvable subalgebra of $\g$ which is stable under $\theta$, since $\theta$ is
involutory. Moreover there exists a Cartan subalgebra $\hh$ of $\g$ which is
contained in
$\ss$ since the intersection of any two Borel subalgebras contains a Cartan
subalgebra. Furthermore from Weyl group theory $$\ss = \hh + \n\cap\n'\eqno
(2.43)$$ is a Levi decomposition of $\ss$. But since $\g(z)$ is stable under
$\theta$ one also has $$\g(z)\s \n\cap\n'\eqno (2.44)$$ But now there exists a
regular semisimple element
$u\in
\hh$ such that the spectrum of $ad\,\u|\n$ is a set of positive numbers. In
particular the spectrum of $ad\,u|\n\cap\n'$ is again strictly positive.
Now let $u' = \theta(u)$ so that $u'\in \hh'$ where $\hh' = \theta(\hh)$.
But since $\ss$ is stable under $\theta$ one has $\hh'\s \ss$.
Interchanging the roles of $\hh$ and $\hh'$ it follows that the spectrum
of $ad\,u'|\n \cap\n'$ is again strictly positive. But, by Lie's theorem, the
adjoint action of $\ss$ on $\n\cap\n$ may be triangularized. The diagonal
entries of both $ad\,u$ and $ad\,u'$ on $\n\cap\n$ are positive. Hence
the same is true of $ad\,v$ where $v= u+ u'$. This however implies that for any
$w\in \n\cap\n'$, 
$$exp\,(-t)\,v\cdot w \,\,\hbox{converges to $0$ as $t$ goes to $+\infty$}\eqno
(2.45)$$ (noting that even though $v$ may not be semisimple the nilpotent
component of $v$ relative to its Jordan decomposition contributes only polynomial
terms in $t$). But this implies that $$\n\cap\n'\s Nil_K\eqno (2.46)$$ since
$v\in \k$. Hence
$z\in
Nil_K$ proving Theorem 1.2. \vs {\bf 2.5.} {\bf Proof of
Theorem 1.3.} That is, we prove that if $z\in \g$ then $z\in Nil_K$ if and only if
$ f(z)=0$ for all homogeneous $f\in A_{2n}$ of positive degree. Of course the ``only
if" is obvious since $A_{2n}\s A$. Assume then that $z\in \g$ and $ f(z)=0$ for all
homogeneous $f\in A_{2n}$ of positive degree. But then recalling the $d(n)\times
d(n)$ matrix $M(z)$ of \S 2.2 one has $$(\xi_z(T_i),\xi_z(T_j))= 0\eqno (2.47)$$
for all $i,j \in \{1,\ldots,d(n)\}$. But then, by Proposition 2.3, one has
$$tr\,ad\,u\,ad\,v= 0\eqno (2.48)$$ for all $u,v\in \g(z)$. Thus, since $ad$
is faithful, $\g(z)$ is solvable and hence its adjoint action on $\g$ can be
triangularized. The nilcone of $\g$ intersected with $\p$ is just the set of zeros
of the polynomials in $S(\p)^K$ of positive degree (see Proposition 11 in [K-R]).
But as one knows the homogeneous generators of $S(\p)^K$ have the same degrees as the
homogeneous generators of the polynomial invariants of the restricted Weyl group
operating on a Cartan subspace of $\p$ (the symmetric space analogue of Chevalley's
theorem). But then one easily has
$S(\p)^K\s A_{2n}$. (This follows, for example, from Proposition 23 in [K-R].) But
since
$S(\k)^K\s A_{2n}$ and
$S(\p)^K\s A_{2n}$ one has that $x$ and $y$ are
nilpotent where
$z=x+y$ is the decomposition (1.13). Thus the diagonal entries of $ad\,x$ and
$ad\,y$ are zero. But since $x$ and $y$ generate $\g(z)$ the diagonal entries of any
element in
$\g(z)$ are zero. Thus any element in $\g(z)$ is nilpotent. Theorem 1.3 then
follows from Theorem 1.2. QED \vs {\bf 2.6.} {\bf Proof of Theorem 1.4.} Theorem 1.3 above and 
Theorem 1.1
in [D] assert that there exists $r$ such that $A=A_r$ where $r\leq max\,\{2, {3\over
8}\,dim\,\p\,(2n)^2\}$. But then Theorem 1.4 follows since ${1\over 2}\,(x(x-1))\geq {3\over 8} \,x^2$ for
$x\geq 4$, and (asuming $\g\neq 0$), one surely has $n>2$. QED\vskip 1.5pc \centerline{\bf Appendix}\vskip
1.5pc The purpose of this appendix is to show that $\g^{K\,reg}$ is
not empty. 

{\bf 1.1A}. Let
$\g = \k +\a +\n$ be a complexified Iwasawa decomposition of $\g$,
consistent with the complexified Cartan decomposition $$\g = \k +
\p\eqno (1.1A)$$ (e.g. $\a$ is a complexified Cartan subspace of
$\p$). Let
$R\s
\a^*$ be the set of restricted roots, and for any $\nu\in R$ let
$\g_{\nu}\s \g$ be the corresponding restricted root space. Let
$R_+\s R$ be the set of positive restricted roots defined
so that $$\n = \oplus_{\nu\in R_+}\,\,\g_{\nu}$$ Let
$\zeta$ be the nonvanishing polynomial function on $\a$ defined by
putting
$$\zeta =\prod_{\nu,\nu'\in R,\,\,\nu\neq
\nu'} (\nu-\nu')\eqno (1.2A)$$ Let $y\in \a$ be defined so that
$$\zeta(y)\neq 0$$ Let $\m$ be the centralizer of $\a$ in $\k$. We
recall that $\theta$ is the complexified Cartan involution
corresponding to (1.1A). For
$\nu\in R_+$ let $x_{\nu}\in \g_{\nu}$. Let $x_{-\nu}\in
\g_{-\nu}$ be defined by putting $x_{-\nu} = \theta\,\,x_{\nu}$.
Let $\widetilde{R} = R\cup\{0\}$ where, here, we regard $0$ as the
zero linear functional on $\a$. Then $\widetilde{R}$
is the set of weights for the adjoint action of $\a$ on $\g$.
Let $\r$ be the $\Bbb C$-span of the set $\{x_{\nu}\}\,\nu \in
\widetilde{R}$. Also let $x = \sum_{\nu\in \widetilde{R}}\,x_{\nu}$
so that $x\in \k$ and also $x\in \r$. \vs {\bf Remark 1.1A.}
Note that, for any $\nu\in R$, $2\,\nu$ is a factor of $\zeta$,  so
that $\nu(y)\neq 0$. \vs Let $z = x +y$ and let $\g(z)$ be the Lie
subalgebra of $\g$ generated by $x$ and $y$. One notes that $\r$ is
stable under $ad\,y$ and that $ad\,y |\r$ is diagonalizable with
distinct eigenvalues. In fact clearly $\r$ is a cyclic $ad\,y$
module with $x$ as cyclic generator and hence \vs {\bf
Proposition 1.2A.} {\it One has $x_{\nu}\in \g(z)$ for any $\nu\in
\widetilde{R}$.} \vs {\bf 1.2A.} The element $y\in \p$ will be
fixed as in \S 1.1A. It will be our objective in this section to show
that $x_0$ and $x_{\nu}, \,\nu\in R_+$ can be chosen, consequently
$x$ can chosen, so that $\g(z) = \g$, i.e. $z\in \g^{K\,reg}$. This
will establish that
$\g^{K\,reg}$ is not empty.  

Let $\hh_{\m}$ be a Cartan subalgebra of $\m$ so that $\hh =
\h_{\m} + \a$ is a Cartan subalgebra of $\g$. Let $\Delta\s
\hh^*$ be the set of roots for $(\hh,\g)$, and for each
$\varphi\in \Delta$, let $e_{\varphi}\in \g$ be a corresponding
root vector. Obviously $\g_{\nu}$ is stable under $ad\,\hh$ for any
$\nu\in R$. Hence there exists a subset $\Delta_{\nu}\s \Delta$
such that $$\g_{\nu} = \sum_{\varphi\in \Delta_{\nu}}\,\,\Bbb
C\,e_{\varphi}\eqno (1.3A)$$ It is immediate that $$\Delta_{-\nu} =
- \Delta_{\nu}\eqno (1.4A)$$ For any $\nu\in R$ let $h_{\nu}\in \a$
be such that, with respect to the Killing form, $(h,h_{\nu}) =
\nu(h)$ for any $h\in \a$. It is clear of course that $\a$ is
spanned by $\{h_{\nu}\mid \nu\in R_+\}$. 

Let $P = {1\over 2}\,(1-\theta)$ so that
$P:\g\to
\p$ is the projection of $\g$ on $\p$ with respect to (1.1A). Since
$\g(z)$ is clearly stable under $\theta$ for any $x\in \k$ it is also
stable under $P$. One easily has \vs {\bf Lemma 1.3A.} {\it Let
$\nu\in R$ and let $\varphi\in \Delta_{\nu}$ so that
$-\varphi\in \Delta_{-\nu}$. Then $$P[e_{\varphi},e_{-\varphi}]
= c\,h_{\nu}\eqno (1.5A)$$ for some $c\in \Bbb C^{\times}$.}\vs
A useful criterion for $K-regularity$ is given in \vs {\bf
Proposition 1.4A.} {\it For
$z$ to be in $\g^{K\,reg}$ it is necessary and sufficient that
$\n\s\g(z)$.} \vs {\bf Proof.} The necessity is by definition.
Assume $\n \s \g(z)$. Then $\g_{\nu}\in g(z)$ for any $\nu\in R_+$.
But clearly $\theta(\g_{\nu}) = \g_{-\nu}$ so that $\g_{-\nu}\s
\g(z)$. But then $h_{\nu}\in \g(z)$ for any $\nu\in R_+$ by Lemma
1.3. Hence $\a +\n\s \g(z)$. But from the Iwasawa decomposition
$P(\a + \n) = \p$. Thus $\p\s \g(z)$. However $\g = \p + [\p,\p]$.
Thus $\g(z) = \g$. QED \vs Let $R_+^1$ be the set of all $\nu\in
R_+$ such that $dim\,g_{\nu} =1$ and let $R_+^2$  be the complement
of $R_+^1$ in $R_+$. Assume $\nu\in R_+^2$. Then the weights of
$ad\,\hh_{\m}$ on $\g_{\nu}$ are of the form $\varphi|\hh_{\m}$
where $\varphi\in \Delta_{\nu}$. Since roots, as weights of
$ad\,\hh$ acting on $\g$, have multiplicity 1 it follows immediately
that the weights of $ad\,\hh_{\m}$ on $\g_{\nu}$ have multipicitity
one. Thus if $\eta_{\nu}$ is the polynomial function	on $\hh_{\m}$
defined by putting $$\eta_{\nu} = \prod_{\varphi,\varphi'\in
\Delta_{\nu},\,\varphi\neq
\varphi'}\,\,(\varphi-\varphi')|\hh_{\m}\eqno (1.6A)$$ then
$\eta_{\nu}$ is nonvanishing. One immediately has \vs
{\bf Proposition 1.5A.} {\it Assume $\nu\in R_+^2$. Let
$x'\in\hh_{\m}$ be such that $\eta_{\nu}(x')\neq 0$. (Such an
element 
$x'$ exists since $\eta_{\nu}$ is nonvanishing.) Then $\g_{\nu}$ is
a cyclic module for $ad\,x'$.}\vs We can now exhibit an
element $z\in \g^{K\,reg}$. Recall the notation of \S1.1A. \vs {\bf
Theorem 1.6A.} {\it For any $\nu\in R_+^1$ let $0\neq x_{\nu}\in 
\g_{\nu}$. If $R_+^2$ is empty let $x_0 = 0$. If $R_+^2$ is not
empty let $\eta$ be the nonvanishing function on $\hh_{\m}$ defined
by putting $$\eta =\prod_{\nu\in R_+^2}\,\eta_{\nu}\eqno (1.7A)$$
Let $x_0\in \hh_{\m}$ be such that $\eta(x_0)\neq 0$ so that (by
Proposition 1.5A) $\g_{\nu}$ is a cyclic module for $ad\,x_0$ for
any $\nu\in R_+^2$. For $\nu\in R_+^2$ let $x_{\nu}\in \g_{\nu}$ be
a cyclic generator of $\g_{\nu}$ with respect to the action of
$ad\,x_0$. Now let $y\in\a$ be as in \S1.1, and as in \S1.1, let $x=
\sum_{\nu\in \widetilde {R}}\,x_{\nu}$ where we recall $x_{-\nu} =
\theta(x_{\nu})$ for $\nu\in R_+$ so that $x\in \k$. Then $\g(z) =
\g$ where $z= x+ y$.} \vs {\bf Proof.} One has
$x_{\nu}\in \g(z)$ for any $\nu\in \widetilde {R}$ by
Proposition 1.2A. Thus $\g_{\nu}\s
\g(z)$ for any $\nu\in R_+^1$. On the other hand if $R_+^2$ is not
empty then
$\g_{\nu}\s \g(z)$ for
$\nu\in R_+^2$ since the Lie algebra generated by $x_0$ and
$x_{\nu}$ contains $\g_{\nu}$. Thus $\n\s \g(z)$ and hence $z\in
\g^{K\,reg}$ by Proposition 1.4A. QED

\vskip 2pc

\centerline{\bf References}\vskip 1.3pc
\parindent=42pt

\item {[D]} H. Derksen, Polynomial bounds for rings of invariants, {\it Proc. Amer.
Math. Soc.}, {\bf 129}, no.4, 955--963
\item {[K-T]} B. Kostant and Juan Tirao, On the structure of certain subalgebras of
a universal enveloping algebra, {\it Trans. Amer. Math. Soc.}, {\bf 218}(1976),
133--154
\item {[K-R]} B. Kostant and S. Rallis, Orbits and Representations associated with
Symmetric Spaces, {\it Amer. J. Math.}, {\bf 93}(1971), No. 3, 753--809 
{\item {[K-W]} H. Kraft and N. Wallach, On the nullcone of representations of
Reductive Groups, {\it Pacific J. Math.}, {\bf 224}(2006), 119--140
\item {[L]} J. Lepowsky, Algebraic results on representations of semisimple Lie
groups, {\it Trans. Amer. Math. Soc.}, {\bf 176}(1973), 1--43
\item {[L-M]} J. Lepowsky and G. McCollum, On the determination of irreducible
modules by restriction to a subalgebra, 
{\it Trans. Amer. Math. Soc.}, {\bf 176}(1973), 44--57
\item {[P-1]} V. Popov, Constructive invariant theory, {\it Ast\'erisque}, {\bf 87-88}(1981), 303--334
\item {[P-2]} V. Popov, The constructive theory of invariants, {\it Math. USSR Izvest.}, {\bf
10}(1982), 359-376
\item {[P-3]} V. Popov, The cone of Hilbert nullforms, {\it Steklov Inst.
Math.} {\bf 241}(2003), 177--194
\item {[V-1]} D. Vogan, The algebraic structure of representations of semi-simple Lie
groups, I, {\it Ann. of Math.,} {\bf 109}(1979), 1--60
{\item {[V-2]} D. Vogan, {\it Representations of Real Reductive Lie Groups,}
Birkh\"auser, PM {\bf 15}(1981)
\item {[W]} N. Wallach, {\it Real Reductive Groups, I}, Academic Press Inc, {\bf
132}, 1988

\smallskip
\parindent=30pt
\baselineskip=14pt
\vskip 1.9pc
\vbox to 60pt{\hbox{Bertram Kostant}
      \hbox{Dept. of Math.}
      \hbox{MIT}
      \hbox{Cambridge, MA 02139}}\vskip 1pc

      \noindent E-mail kostant@math.mit.edu

\end

\end

\end